\input amstex
\documentstyle{amsppt}
\magnification=1200
\pagewidth{16.8truecm}
\pageheight{23truecm}
\nologo

\topmatter

\title Timelike surfaces with 
harmonic inverse mean curvature 
\endtitle
\leftheadtext{Fujioka and Inoguchi}
\rightheadtext{Harmonc inverse mean curvature surfaces}
\author Atsushi Fujioka
\footnote"\dag"
{Partially supported by Grant-in-Aid for
Encouragement of Young Scientists
No.~12740037,
\newline
\leftline{ Japan Society for Promotion
of Science}
}
and Jun-ichi Inoguchi
\footnote"*"
{Partially supported by Grant-in-Aid for
Encouragement of Young Scientists
No.~12740051,
\newline 
\leftline{Japan Society for Promotion
of Science}
\leftline
\newline
\leftline
{to appear in: Proceeding of the 9th MSJ-IRI,
{\it Integrable Systems in Differential Geometry, Tokyo, 2000.}} 
 }
\endauthor
\affil Kanazawa University and Fukuoka University
\endaffil
\address 
(Fujioka) 
Department of Mathematics,
Faculty of Science,
Kanazawa University,
\flushpar
Kakuma-machi,
Kanazawa,
Ishikawa,
920--1192,
Japan
\endaddress
\email fujioka\@kenroku.kanazawa-u.ac.jp
\endemail
\address
Current address:
Graduate School of Economics,
Hitotsubashi University, 2-1,
Naka, Kunitachi, Tokyo,
186--8601,
Japan
\endaddress
\email fujioka\@math.hit-u.ac.jp
\endemail

\address
(Inoguchi)
Deparment of Applied Mathematics,
Fukuoka University,
Fukuoka, 
\flushpar
814-0180,
Japan \endaddress 
\email
inoguchi\@bach.sm.fukuoka-u.ac.jp
\endemail
\address
Current address: 
Department of Mathematics Education,
Utsunomiya University, 
Utsunomiya, 321-8505,
Japan
\endaddress
\email
inoguchi\@cc.ustunomiya-u.ac.jp
\endemail

%
\endtopmatter
\document

\head Introduction
\endhead

In this paper we introduce the 
notion of {\it timelike surface
with harmonic inverse mean curvature}
in $3$-dimensional 
Lorentzian space forms,
and study their fundamental properties.
\smallskip
\par

In classical differential geometry,
surfaces of constant mean curvature
(CMC surfaces) have
been studied extensively \cite{1}.
As a generalization of CMC surfaces,
Bobenko \cite{2} introduced the
notion of 
{\it surface with harmonic inverse 
mean curvature} (HIMC surface).
He showed that
HIMC surfaces admit 
Lax representation
with {\it variable} spectral
parameter.
In \cite{5}, Bobenko, Eitner and Kitaev
showed that
the Gauss equations of $\theta$-isothermic HIMC
surfaces reduce to the ordinary differential equation:
$$
\left(
\frac{q^{\prime \prime}(t)}
{q^{\prime}(t)}
\right )^{\prime}
-q^{\prime}(t)=
{\Cal S}(t)
\left (
2-\frac{q^{2}(t)+c}{q^{\prime}(t)}
\right),\ \
q^{\prime}(t)<0,
\tag$*$
$$ 
with $c=\theta^2>0$. 
Here the coefficient function 
${\Cal S}(t)$ is 
$1/\sin^{2}(2t)$,
$1/\sinh^{2}(2t)$
or $1/t^{2}$.
This ordinary differential equation is called the 
{\it generalized Hazzidakis equation}.
Bobenko, Eitner and Kitaev
\cite{5}
solved $(*)$ in terms of
Painlev{\'e} transcendents $P_{\roman V}$
and $P_{{\roman V}\!{\roman I}}$.

For $c<0$, solutions to $(*)$
do not describe surfaces in Euclidean $3$-space.
It seems to be interesting to find 
``corresponding surfaces" to such solutions.

\smallskip
\par
The first
author extended the notion of
HIMC surface in Euclidean 
$3$-space to that of 
Riemannian 3-space forms \cite{7}.
Moreover he generalized a theorem due to
Lawson ({\it Lawson correspondence})
to HIMC surfaces.
By using the Lawson correspondence
for HIMC surfaces,
we have classified Bonnet surfaces 
with constant curvature in 
Riemannian 3-space forms
\cite{8}.
Corresponding results for spacelike
surfaces in Lorentzian 3-space forms are
obtained in \cite{10}.

\smallskip
\par
On the contrary very little
is known about (integrable)
timelike surfaces  of 
nonconstant mean curvature
in Lorentzian 3-space forms.
Timelike Bonnet surfaces are investigated by 
present authors very recently \cite{11}.

\smallskip
\par
In this paper we introduce the notion of
timelike surface with harmonic inverse
mean curvature
(THIMC surface)
in Lorentzian
3-space forms.
We shall show that every 
solution to the generalized Hazzidakis
equation with $c<0$ describes
a THIMC surface in Minkowski
3-space.
This is one of the motivations to
study THIMC surfaces. 

Because of the indefiniteness of metric, 
timelike surface geometry has many 
aspects different from 
Euclidean surface geometry. 
For instance, there exist 
timelike (HIMC) surfaces with 
imaginary principal curvatures.
Moreover there exist non totally umbilical
timelike surfaces with
real repeated principal curvatures. 
Both of such surfaces 
have no counterparts 
in Euclidean surface geometry
and spacelike surface geometry.
Thus the geometry of THIMC
surfaces has its own interest.

\smallskip
\par
The second motivation of the present study is
to give new examples
of Lax equations with variable spectral 
parameter,
namely,  Lax equations 
whose spectral parameters depend on the
variables.
Burtsev, Zakharov and Mikhailov \cite{6}
exhibited some examples of Lax equations
with variable spectral parameter appeared
in theoretical physics.
In differential geometry, HIMC surfaces and Bianchi
surfaces are known examples. 
(See \cite{2}, \cite{15}
and  \cite{16}).

We shall show that
THIMC surfaces
 Lorentzian 3-space forms
admit Lax representation 
with variable spectral
parameter. Moreover we shall show that
in de Sitter 3-space or anti de
Sitter 3-space, THIMC surfaces
admit Lax representation with {\it two} 
independent variable spectral
parameters. 

\smallskip
\par
This paper is organized as follows.
After recalling fundamental
facts on Lorentzian geometry,
we introduce the notion of THIMC surface
in Minkowski 3-space in Section 3.
We give a Lax representation and 
an immersion formula (Sym-formula)
for THIMC surfaces.
Some elementary examples
will be given in Section 3.
In the next Section 4,
we introduce the notion of
$\pm$ isothermic timelike surface.
We shall give a duality between
timelike Bonnet surfaces and  
$\pm$ isothermic 
THIMC surfaces.
\par
In Section 5, we shall investigate
the normal forms of 
the Gauss equations of THIMC surfaces.
More precisely we show that 
($\theta$-isothermic or anti $\theta$-isothermic)
THIMC surfaces in Minkowski 3-space are derived from solutions
to the generalized Hazzidakis equation with $c=-\theta^2<0$.
 
In Section 6, we shall generalize
the notion of THIMC surface
to Lorentzian 3-space forms
and establish a Lawson-type  correspondence
for THIMC surfaces.
 
\smallskip

The authors would like to
express their gratitude
to the referee for 
careful reading of the
manuscript.

\head 1. Lorentzian space forms
\endhead
\subhead 1.1
\endsubhead
First of all, we shall describe 
{\it Lorentzian $3$-space forms},
i.e., 
complete and
connected
\flushpar
 Lorentzian 
$3$-manifolds
${\frak M}^3_1(c)$ of
constant curvature $c$ 
explicitly.
\par
Without loss of generality, we may assume that
$c=0$ or $\pm 1$.

On a Cartesian 4-space
${\bold R}^4$,
we equip the following scalar product
$\langle \cdot,\cdot \rangle=
\langle \cdot,\cdot \rangle_c$:

$$
\align
& \langle a,b \rangle_c=-a_0b_0+a_1b_1+a_2b_2
+a_3b_3, \ \ c=1,  \\
& \langle a,b \rangle_c=-a_0b_0-a_1b_1+a_2b_2
+a_3b_3, \ \ c=0, \ -1.
\endalign
$$
The resulting semi-Euclidean 4-space
$({\bold R}^4, \langle\cdot,\cdot\rangle)$
is of index $1$ for $c=1$ and of index $2$
for $c=0$ or $-1$ respectively. 
The Lorentzian 3-space
forms ${\frak M}^3_1(c)$ are embedded in
the semi-Euclidean space 
$({\bold R}^4,\langle \cdot, \cdot \rangle_c)$ 
as
$$
\align
& {\frak M}^3_1(0)=\{p\in
({\bold R}^4,\langle \cdot, \cdot \rangle_0)
\ \vert \ p_0=0
\}={\bold E}^3_1, \ \ 
\text{the Minkowski 3-space}, \\
& {\frak M}^3_1(1)=\{p\in
({\bold R}^4,\langle \cdot, \cdot \rangle_1)
\ \vert \ \langle p,p \rangle_1=1 \
\}=S^3_1, \ \ \text{the de Sitter
3-space}, \\
& {\frak M}^3_1(-1)=\{p\in
({\bold R}^4,\langle \cdot, \cdot \rangle_{-1})
\ \vert \ \langle p,p \rangle_1=-1 \
\}=H^3_1, \ \ 
\text{the anti de Sitter 3-space} 
\endalign
$$

For more details on semi-Riemannian geometry, we refer to 
O'Neill \cite{18}.

\subhead 1.2
\endsubhead
Next we recall 2 by 2 matrix models of
${\frak M}^3_1(c)$ for later use.

First the semi-Euclidean 4-space
${\bold E}^4_2=({\bold R}^4,\langle \cdot,\cdot \rangle_{-1})$
is identified with the linear space ${\roman M}_2{\bold R}$ 
of all 2 by 2 real matrices via the isomorphism:

$$
p=(p_0,p_1,p_2,p_3)
\longleftrightarrow
p_0{\bold 1}+p_1{\bold i}+p_2{\bold j}'+p_3{\bold k}'
=\pmatrix
p_0-p_3 & -p_1+p_2 \\
p_1+p_2 & p_0+p_3
\endpmatrix. \tag{1.1}
$$

The semi-Euclidean metric of 
${\bold E}^4_2$ corresponds to the following scalar product
on ${\roman M}_{2} {\bold R}$.
$$
\langle X,Y \rangle=\frac{1}{2}\{ {\roman t}{\roman r}(XY)
-{\roman t}{\roman r}(X){\roman t}{\roman r}(Y)\},
\ \ 
X,Y \in {\roman M}_2 {\bold R}.
\tag{1.2}
$$
Under the identification (1.1),
the Minkowski 3-space ${\bold E}^3_1(p_1,p_2,p_3)$
is identified with
the Lie algebra 
${\frak g}={\frak s}{\frak l}_{2}{\bold R}$:
$$
{\frak g}= \{X \in {\roman M}_{2}{\bold R}\
 | \
\roman{tr}\ X=0 \
 \}
$$
with metric 
$$
\langle X,Y \rangle=\frac{1}{2}{\roman t}{\roman r}
(XY),\ X,Y \in {\frak g}.
$$
\subhead 1.3
\endsubhead
Next, since  
$$
\langle X, X \rangle=-\det X
$$
for all $X \in {\roman M}_{2}{\bold R}$,
the anti de Sitter 3-space $H^3_1 \subset {\bold E}^4_2$
corresponds to the real special linear group:
$$
G={\roman S}{\roman L}_2{\bold R}=\left \{ 
\pmatrix
a & b \\
c & d
\endpmatrix \in {\roman M}_{2}{\bold R}
\ \biggr | \ ad-bc=1 \right \}. 
$$

Since the Lorentzian metric of $G$ is 
bi-invariant,
the product group $G \times G$ acts transitively and
isometrically on 
$H^3_1$ as follows:

$$
\mu_{H}:(G \times G )\times H^3_1 \longrightarrow  H^3_1,\ \ 
\mu_{H}(g_1,g_2)X=g_1\ X \ g_2^{-1}
$$
for $(g_1,g_2)\in G \times G,\ X \in H^3_1$.
The isotropy subgroup $\Delta$ of $G \times G$ at ${\bold 1}$
is the diagonal subgroup of $G \times G$, that is,
$\Delta=\{(g_1,g_1 ) \  | g_1 \in G \}$.
Hence the anti de Sitter $3$-space $H^3_1$ is represented
by $H^3_1=(G \times G)/ \Delta$
as a Lorentzian symmetric space.
The natural projection $p_H: G\times G \rightarrow H^3_1$
is given explicitly by
$p_H(g_1,g_2)=g_1\
g_2^{-1},\ (g_1,g_2) \in G \times G.$
\bigskip
\par
Moreover $G$ acts isometrically on ${\bold E}^3_1$ via the
Ad-action:
$$
\roman {Ad}:G \times {\bold E}^3_1 \to 
{\bold E}^3_1;\ \roman{Ad}(a)X=aXa^{-1},\
a \in G, \ X \in {\bold E}^3_1.
$$
\medskip
\par
\subhead 1.4
\endsubhead
Finally
we recall a 2 by 2 matrix
model
of $S^3_1$.
The Minkowski 4-space
${\bold E}^4_1=({\bold R}^4,\langle \cdot,\cdot
\rangle_{1})$ is identified with
the space
${\Bbb H}$ of all Hermitian $2$-matrices
via the following isomorphism:
$$
p=(p_0,p_1,p_2,p_3)
\longleftrightarrow
\pmatrix
p_0+p_1 & p_3-\sqrt{-1}p_2  \\
p_3+\sqrt{-1}p_2 & p_0-p_1 
\endpmatrix 
\in {\Bbb H}.
\tag1.4
$$

Under the identification (1.4),
the scalar product $\langle\cdot,\cdot 
\rangle_{1}$
of ${\bold E^4_1}$
corresponds to
the following scalar product
on ${\Bbb H}$:
$$
\langle X,Y \rangle
=-\frac{1}{2}{\roman t}{\roman r}
({\bold i}^{\prime}X\ {\bold i}^{\prime}Y^{t}),\ X,Y \in {\Bbb H},
\ \ 
{\bold i}^{\prime}
=\pmatrix 
0 & -\sqrt{-1} \\
\sqrt{-1} & 0
\endpmatrix.
\tag{1.5}
$$
In particular 
$\det X=-\langle X,X \rangle_{1}$ under (1.4).
Thus the  de Sitter $3$-space $S^3_1$ is represented by
$$
S^3_1=\{ X \in \Bbb H \ |\ \det X=-1 \}.  
$$
The complex special linear group ${\roman S}{\roman L}_2{\bold C}$
acts transitively and isometrically on $S^3_1$ by
$$
\mu_{S}: {\roman S}{\roman L}_2{\bold C} \times
S^3_1 \rightarrow S^3_1,\ \ 
\mu_{S}(g)=gXg^{*}.
$$
Here $g^{*}$ denotes the transposed complex conjugate 
of $g$. 
The isotropy subgroup of 
${\roman S}{\roman L}_2{\bold C}$
at ${\bold i}^\prime$ is
${\roman S}{\roman L}_2{\bold R}$.
Hence the de Sitter $3$-space
$S^3_1$ is represented by
$S^3_1=G^{\bold C}/G$
as a Lorentzian symmetric space.
The natural projection $p_{S}:G^{{\bold C}} \rightarrow 
G^{{\bold C}}/G$ is given explicitly by
$p_{S}(g)=\mu(g)\ {\bold i}'=g\ {\bold i}'\ g^{*},\ g \in G^{{\bold C}}$.

\head 2. Timelike surfaces in Lorentzian  space forms.
\endhead 
We start with some preliminaries on 
the geometry of timelike surfaces
in 
Lorentzian space forms ${\frak M}^3_1(c)$. 

\subhead 2.1
\endsubhead
Let $M$ be a connected $2$-manifold 
and $F: M \rightarrow {\frak M}^3_1(c)$
an immersion. The immersion $F$ is said to be {\it timelike} if
the induced metric $I$ of $M$ is Lorentzian.
Hereafter we may assume that $M$ is an orientable timelike surface in
${\frak M}^3_1(c)$ immersed by $F$.
The induced Lorentzian metric $I$ of $M$
determines a Lorentzian conformal
structure on $M$.
We treat $M$ as a Lorentz surface with respect to
this conformal structure
and $F$ as a conformal immersion. 
Our general reference on ``Lorentz surfaces" is Weinstein \cite{21}.

On a timelike surface $M$, there exists a local coordinate system
$(x,y)$ such that
$$
I=e^{\omega}(-dx^2+dy^2). 
\tag{2.1}
$$
Such a local coordinate system $(x,y)$ is called a
{\it Lorentz isothermal coordinate system}.

Let $(u,v)$ be the
local {\it null coordinate system}
of $M$
derived from $(x,y)$. Namely $(u,v)$ is defined by
$u=x+y,\ v=-x+y$.
Then the induced metric $I$ can be written as 
$$
I=e^{\omega}dudv. \tag{2.2}
$$
Now, let $N$ be a unit normal vector field to $M$.
The {\it second fundamental form} $I\!I$ of $(M,F)$ 
derived from $N$ is defined by
$$
I\!I=-\langle dF,dN \rangle.
$$
The {\it shape operator}
$S$ of $(M,F)$ relative to $N$ is
defined by
$$
S=-dN.
$$
The (complex) eigenvalues of
$S$ are called {\it principal curvatures} of
$(M,F)$.
The {\it mean curvature} $H$ of
$(M,F)$ is defined by $H=\roman{tr}\ S/2$. 
The Gaussian curvature $K$ of $(M,I)$ is computed by
the formula: $K=\det S$.

\smallskip
\par
The
Gauss-Codazzi equations
of $(M,F)$
have the following form:
$$
\omega_{uv}+\frac{1}{2}(H^{2}+c)e^{\omega}
-2QRe^{-\omega}=0, \tag{G$_{c}$}
$$
$$
H_{u}=2e^{-\omega}Q_v,\ \ 
H_{v}=2e^{-\omega}R_u.\tag{C$_{c}$}
$$
Here
the functions $Q=\langle F_{uu},N \rangle$
and  $R=\langle F_{vv},N \rangle$
define global 
null $2$-differentials 
$Q^{\#}=Qdu^2$ and $R^{\#}=Rdv^2$
on $M$. These two differentials 
are called the {\it Hopf differentials}
of $M$.
The Gauss equation implies 
$$
K=-2\ \omega_{uv}\ e^{-\omega}.
\tag2.3
$$
Let us denote by
${\Cal D}$ the discriminant
of the characteristic equation:
$$
\det (t{\roman I}-S)=0
$$
for the shape operator $S$.
Here ${\roman I}$ is the identity
transformation of the tangent
bundle $TM$ of $M$.
Then by the Gauss equation, we have 
$$
{\Cal D}=H^2-K+c=4e^{-2 \omega}\ Q R.
\tag2.4
$$
The first and second fundamental forms 
are related by the formula:
$$
I\!I-HI=Q^{\#}+R^{\#}.
$$
This formula implies that
the common zero of $Q$ and $R$ coincides with the
umbilic point of $(M,F)$.
Even if $S$ has real and same eigenvalues,
$(M,F)$ is not necessarily totally umbilic.
In fact, there exist timelike surfaces with $QR=0$ but $I\!I-HI\not=0$.
See Example 3.3.

\subhead 2.2
\endsubhead
In the study of timelike surfaces, we also use 
the following 
local coordinate system:
\proclaim{Lemma 2.1}
Let $F:M \rightarrow {\frak M}^3_1(c)$
be a timelike surface.
Then there exists a local coordinate system $({\check u},{\check v})$
such that 
$$
I=-e^{\check \omega}d{\check u}d{\check v}.
\tag 2.5
$$
With respect to this coordinate system, the Gauss-Codazzi
equations are written as
$$
{\check \omega}_{{\check u}{\check v}}
-\frac{1}{2}(H^{2}+c)e^{\check \omega}
+2{\check Q}{\check R}e^{-{\check \omega}}=0, 
\tag{G$_{c}^{-}$}
$$
$$
H_{\check u}=-2e^{-{\check \omega}}{\check Q}_{\check v},\ \ 
H_{\check v}=-2e^{-{\check \omega}}{\check R}_{\check u}
\tag{C$_{c}^{-}$}
$$
for 
${\check Q}=\langle F_{{\check u}{\check u}},N \rangle,
{\check R}=\langle F_{{\check v}{\check v}},N \rangle$.
\endproclaim
We call the local coordinate system $({\check u},{\check v})$ an
{\it anti isothermal coordinate system}.
Anti isothermal coordinate systems will be used for
introducing the
notion of the Christoffel transformation of 
an anti isothermic surface. See Proposition 4.14.

\head 3. Timelike HIMC surfaces in Minkowski $3$-space 
\endhead

In this section we shall consider a generalization
of timelike CMC surfaces in
\flushpar
Minkowski $3$-space
in terms of 
integrability theory.
\subhead 3.1
\endsubhead
We start with recalling 
the Lax representation for
timelike surfaces in ${\bold E}^3_1$.
Hereater we assume $H\not=0$.

Let $F: M \rightarrow {\bold E}^3_1$ be a timelike
surface.
Let us take an SL$_{2}{\bold R}$-valued framing
$\Phi$ defined by
$$
{\roman A}{\roman d}(\Phi)({\bold i},\ {\bold j}',\ {\bold k}')=
(e^{-\frac{\omega}{2}}F_x,
\  e^{-\frac{\omega}{2}}F_y,\ N).
$$
Thus we get the following Lax representation
of Gauss-Codazzi equations:

$$
\frac{\partial}{\partial u}\Phi=\Phi U,\ 
\frac{\partial}{\partial v}\Phi=\Phi V,
\tag{3.1}
$$
$$
U=\pmatrix
-\frac{1}{4}\omega_{u} & -Qe^{-\frac{\omega}{2}} \\
\frac{H}{2}e^{\frac{\omega}{2}} & \frac{1}{4}\omega_u
\endpmatrix,\ \ 
V=\pmatrix
\frac{1}{4}\omega_v & -\frac{H}{2}e^{\frac{\omega}{2}} \\
Re^{-\frac{\omega}{2}} & -\frac{1}{4}\omega_v
\endpmatrix. \tag{3.2}
$$

Now we shall insert a {\it variable spectral
parameter} $\lambda$ , {\it i.e.}, an additional
{\it real \/} parameter $\lambda$ depends on the 
coordinate system
$(u,v)$ into the Lax pair (3.2) 
in the following way:
$$
U_{\lambda}=\pmatrix
-\frac{1}{4}\omega_u & -Qe^{-\frac{\omega}{2}} \\
\frac{H}{2}{\lambda}e^{\frac{\omega}{2}} & \frac{1}{4}\omega_u 
\endpmatrix, \
V_{\lambda}=\pmatrix
\frac{1}{4}\omega_{v} & -\frac{H}{2}\lambda^{-1} e^{\frac{\omega}{2}} \\
Re^{-\frac{\omega}{2}} & -\frac{1}{4}\omega_v
\endpmatrix.\tag{3.3}
$$ 
 
Then the compatibility condition 
$$
\frac{\partial}{\partial u}V_{\lambda}
-\frac{\partial}{\partial v}U_{\lambda}
+[\ U_{\lambda},V_{\lambda}\ ]=0 \tag{3.4}
$$
for the deformed Lax pair
$\{U_{\lambda},V_{\lambda}\}$ yields

$$
\omega_{uv}+\frac{1}{2}H^{2}e^{\omega}-2QRe^{-\omega}=0, 
\tag{${\roman G}_0$}
$$

$$
Q_v=\frac{e^{\omega}}{2}(H\lambda^{-1})_{u},\ \ 
R_u=\frac{e^{\omega}}{2}(H\lambda)_{v}.
\tag{3.5}
$$
The Lax pair 
$\{U_{\lambda},V_{\lambda}\}$
describes a timelike surface in ${\bold E}^3_1$
if and only if the equations (3.5) are consistent
with Codazzi equations (${\roman C}_0$).
The equations (3.5) is consistent with
(${\roman C}_0$) if and only if
$$
\frac{\partial}{\partial v}\{ H(1-\lambda)\}=0,\ \
\frac{\partial}{\partial u}\{ H(1-\lambda^{-1})\}=0.
\tag{3.6}
$$
These equations (3.6) can be easily solved as follows:
$$
H=\frac{1}{f(u)+g(v)},\ \ \ \lambda=-\frac{g(v)}{f(u)},
\tag{3.7}
$$
where $f(u)$ and $g(v)$ are smooth functions.
It is easy to see that the mean curvature $H$ is invariant
under the one parametric deformation

$$
f \longmapsto f+\frac{1}{2\tau},\ \ g 
\longmapsto g-\frac{1}{2\tau}
,\ \ \ \tau \in {\bold R^{*}}.
$$

Under this deformation, the spectral parameter $\lambda$
is transformed as
$$
\lambda=-\frac{g}{f}
\longmapsto \lambda(u,v;\tau)=\frac{1-2\tau g}  
{1+2\tau f},\ \ \ \tau \in  {\bold R}.
$$
Note that
$\lambda(u,v;0)\equiv 1$.
The form (3.7) of $H$ is equivalent to
the Lorentz-harmonicity of $1/H$,
{\it i.e.,\/} $(1/H)_{uv}=0$.
As in the Euclidean surface geometry \cite{2}
and spacelike surface geometry \cite{10},
we shall call a timelike surface $M$ in 
${\bold E}^3_1$, a {\it timelike
surface with harmonic inverse mean curvature}
(THIMC surface)
if $1/H$ is a Lorentz-harmonic function.

\bigskip
\subhead 3.2
\endsubhead
Here we would like to exhibit three
elementary examples of THIMC surfaces.
\medskip
\par
\example{Example 3.1} (THIMC cylinders.)
Let $a(y)=(a_2(y),a_3(y))$ be a curve in Euclidean
plane ${\bold E}^2(\xi_2,\xi_3)$ parametrized
by the arclength parameter $y \in {\Cal I}$. Here ${\Cal I}$ 
is an interval. {\it A timelike cylinder over the curve}
$a$ is a flat timelike surface in ${\bold E}^3_1$
defined by the immersion
$F:{\Cal I} \times {\bold R} \longrightarrow
{\bold E}^3_1;\ \ F(x,y)=
(x,a_2(y),a_3(y))$.
It is straightforward to see that the mean curvature
of the cylinder is
$H=\kappa(y)/2$. 
Here 
$\kappa$ is the curvature of 
$a$. 
Thus the cylinder $F$ is a THIMC surface if and
only if the base curve has the curvature
$\frac{1}{C_1 y+C_2}$, $C_1, C_2 \in {\bold R}$.
It is well known that curves with curvature 
$\frac{1}{C_1 y+C_2}$
are logarithmic spirals or circles.
Hence all the THIMC cylinders over a Euclidean
curve are cylinders over a
logarithmic spiral or a circular cylinder. 
\endexample

\example{Example 3.2} (THIMC cylinders over timelike curves.)
Let $a(x)=(a_1(x),a_2(x))$ be a timelike curve in Minkowski
plane ${\bold E}^2_1(\xi_1,\xi_2)$ pa\-ram\-e\-trized
by the proper time parameter $x$ defined on an
interval ${\Cal I}$. 
{\it A timelike cylinder over the timelike curve}
$a$ is a flat timelike surface in 
${\bold E}^3_1$
defined by the immersion
$F:{\Cal I} \times {\bold R} \longrightarrow
{\bold E}^3_1;\ \ F(x,y)=
(a_1(x),a_2(x),y)$.
The mean curvature of $F$ is
$H=\kappa(x)/2$.
Here $\kappa$ is the curvature of $a$.
The cylinder $F$ is THIMC if and only if
$1/\kappa=C_1 x+C_2$, $C_1, C_2 \in {\bold R}$.

We can see that
timelike curves with curvature $\frac{1}{C_1 x+C_2}$
are logarithmic pseudo-spirals or timelike hyperbolas.
({\it cf.\/} Appendix of \cite{10}.)
Hence all the THIMC cylinders are cylinders over a
logarithmic pseudo-spiral or a timelike hyperbola. 
\endexample

\example{Example 3.3}($B$-scrolls.)
A curve $\gamma(s)$ in 
${\bold E}^3_1$
is said to be a {\it null Frenet curve\/}
if it admits a frame field
${\Cal L}=(A,B,C)$
along $\gamma$ 
(called a {\it null frame field\/})
such that $A=\gamma^{\prime}$,
$$
\langle A,A \rangle=
\langle B,B \rangle=0, \ \ 
\langle A,B \rangle=1,\ \ 
\langle C,C \rangle=1, \ \ 
\langle A,C \rangle=
\langle B,C \rangle=0,
$$
$$
\frac{d}{ds}{\Cal L}=
{\Cal L}
\pmatrix
0 & 0 & -\tau \\
0 & 0 & -\kappa \\
\kappa & \tau & 0
\endpmatrix.
$$
The functions $\kappa$ and $\tau$
are called the {\it curvature} and {\it torsion}
of $\gamma$ respectively.
The ruled surface 
$F(s,t)=\gamma(s)+tB(s)$
is called the $B$-{\it scroll}
of $\gamma$.
(See Graves \cite{12} and McNertney \cite{17}).
The mean curvature of $F$ is
the torsion $\tau(s)$. 
It is straightforward to check that
for any
null Frenet curve with $\tau\not=0$,
its $B$-scroll is a THIMC surface.
\endexample

\example{Remark}
The Gaussian curvature
of the $B$-scroll is $\tau^2$.
Thus every $B$-scroll satisfies $QR=0$ but 
is not totally umbilical
($I\!I-HI\not=0$).
The property $QR=0$ implies that every $B$-scroll 
is a timelike Bonnet surface.
Here {\it timelike Bonnet surfaces} 
are timelike surfaces which 
admit nontrivial isometric deformation preserving
mean curvature \cite{11}.

Conversely we proved that 
every timelike Bonnet surface with 
$QR=0$ are $B$-scrolls
\cite{11}. 
\endexample
\smallskip

\subhead 3.3
\endsubhead
In \cite{14}, we have obtained a 
one-parameter ``isometric" deformation of
timelike surfaces with constant mean curvature
(TCMC surfaces). For THIMC surfaces in 
${\bold E}^3_1$,
we get the following 
one-parameter family of ``conformal" deformation.

\proclaim{Proposition 3.4}
Let $F: M \rightarrow {\bold E}^3_1$ be a timelike
surface with harmonic inverse mean curvature.
Express the mean curvature $H$ as
$$
H=\frac{1}{f(u)+g(v)}
$$
in terms of null coordinate
system $(u,v)$.
Here $f(u)$ and $g(v)$ are smooth functions.
Then $F$ admits the following Lax
representation with variable
spectral parameter $\lambda(u,v;\tau)
=(1-2\tau g(v))/(1+2\tau f(u)),\ \tau 
\in {\bold R}
${\rm :}

$$
\frac{\partial}{\partial u}\Phi_{\lambda}=\Phi_{\lambda}U_\lambda,\ 
\frac{\partial}{\partial v}\Phi_{\lambda}=\Phi_{\lambda}V_\lambda,
\tag{3.8}
$$

$$
U_{\lambda}=\pmatrix
-\frac{1}{4}\omega_{u} & -Qe^{-\frac{\omega}{2}} \\
\frac{H}{2}\lambda e^{\frac{\omega}{2}} & \frac{1}{4}\omega_{u} 
\endpmatrix,\ \ 
V_{\lambda}=\pmatrix
\frac{1}{4}\omega_{v} & -\frac{H}{2}\lambda^{-1} e^{\frac{\omega}{2}} \\
R e^{-\frac{\omega}{2}} & -\frac{1}{4}\omega_v
\endpmatrix.
$$

Let $\Phi_{\lambda}(u,v)$ be a solution of {\rm (3.8)}. Then

$$
F_{\lambda}=-
\frac{\partial}{\partial \tau}\Phi_{\lambda}\cdot 
\Phi_{\lambda}^{-1}
\tag{3.9}
$$
describes a family of THIMC
surfaces through $F=F_\lambda \vert_{\tau=0}$ with Gauss map
$N_{\lambda}={\roman A}{\roman d}
(\Phi_{\lambda})\ {\bold k}^{\prime}$.
The fundamental associated quantities 
of $F_{\lambda}$ are given as follows{\rm :}

$$
I_{\lambda}=\frac{e^{\omega}dudv}
{(1+2\tau f)^2(1-2\tau g)^2}, \tag{3.10} 
$$

$$
\frac{1}{H_\lambda}=f_{\lambda}+g_{\lambda}, \ 
f_{\lambda}=\frac{f}{(1+2\tau f)},\ 
g_{\lambda}=\frac{g}{(1-2\tau g)},\ \
\tag{3.11}
$$
$$
Q_{\lambda}=\frac{Q}{(1+2\tau f)^2},\ \ 
R_{\lambda}=\frac{R}{(1-2\tau g)^2},\tag{3.12}
$$

$$
K_{\lambda}=(1+2\tau f)(1-2\tau g)K,\tag{3.13} 
$$

$$
H_{\lambda}^2/K_{\lambda} \equiv H^2/K.
\tag{3.14}
$$

\endproclaim

The formula (3.14) implies that the members of 
the one parameter family
$F_{\lambda}$ have the same ratio of the principal curvatures.


\head 4. $\pm$ isothermic timelike surfaces
\endhead
\subhead 4.1
\endsubhead
In the study of HIMC surfaces in Riemannian space
forms, {\it isothermic surfaces} play a fundamental role.
In this section we shall consider such surfaces 
in timelike surface geometry.

\example{Definition 4.1}
Let $F
: M \rightarrow {\frak M}^3_1(c)$ be a timelike
surface. Then $(M,F)$ is said to be {\it isothermic} 
if there exists a local 
isothermal--curvature line coordinate
system
around any point of $M$.
\endexample

Here an isothermal--curvature line coordinate system is a
local Lorentz--isothermal coordinate 
system such that both of parameter
curves are curvature lines. 
It should be remarked that isothermic property implies
the positivity of the descriminant ${\Cal D}$
for the characteristic equation for the
shape operator $S$.
Equivalently, every isothermic timelike surface 
has real distinct principal curvatures.

The isothermic property for timelike surfaces in 
${\frak M}^3_{1}(c)$ can be reformulated in terms
of {\it associated null coordinate
system} as follows.

\proclaim{Proposition 4.2}
A timelike surface $(M,F)$
is isothermic
if and only if there exists a local null coordinate
system $(u,v)$ 
around any point of $M$
such that 
the Hopf differentials 
take the following form{\rm :}
$$
Q(u,v)=\frac{1}{2}{\frak q}(u,v)\varrho(u)
,\ \ 
R(u,v)=\frac{1}{2}{\frak q}(u,v)\sigma(v),
\ \ \varrho>0,\ \sigma>0.
\tag{4.1}
$$ 
Here ${\frak q}$ is a real smooth function and 
$\varrho$ and $\sigma$ are positive
Lorentz holomorphic and anti holomorphic
functions respectively.
\endproclaim

\example{Remark} On a Lorentz surface $M$ with null coordinate system 
$(u,v)$, a smooth function $f$ on $M$ depends only on $u$
[resp\. $v$] is called a  {\it Lorentz holomorphic function}
[resp\. {\it Lorentz anti holomorphic function}].
\endexample
Hereafter we shall call a null coordinate system
derived from an isothermic coordinate system simply
an {\it isothermic coordinate system}. 

\example{Remark}
Isothermic timelike 
surfaces in ${\bold E}^3_1$ correspond to
solutions of the Zoomeron equation 
studied in
soliton theory.
Note that Zoomeron equation is 
related to Davey-Stewartson 
$\roman{I}\!\roman{I}\!\roman{I}$-equation.
See  Schief \cite{19, p.~97}. 
\endexample

\medskip
\par
\subhead 4.2
\endsubhead
Typical examples of isothermic timelike surfaces are 
timelike surfaces of revolution in ${\bold E}^3_1$.
Here we recall the notion of timelike surfaces
of revolution in
${\bold E}^3_1$. 
A {\it revolution of \/} ${\bold E}^3_1$
is a linear isometry which lies in the identity component
${\roman O}^{++}_1(3)$ of the Lorentz group O$_1(3)$. 
Every revolution
fixes a line pointwise. Such fixed line of a revolution is called
the {\it axis of \/} revolution. Hence revolutions of
${\bold E}^3_1$ can be characterised by the causal character of
the axis.

By a {\it timelike surface of revolution} in ${\bold E}^3_1$
we mean a timelike surface obtained by
revolving about an axis a 
regular curve lying
in some plane containing the axis
\cite{17}.

\example{Example 4.3}({\bf Spacelike axis and Euclidean profile curve.})
Let $F:M \longrightarrow {\bold E}^3_1$ be a timelike
surface of revolution with spacelike axis
and Euclidean profile curve.
Then there exists an isothermic parametrization 
$$
F(x,y)=\frac{1}{a}\left(
e^{\frac{\omega(y)}{2}}\sinh (ax),
e^{\frac{\omega(y)}{2}}\cosh (ax),
c(y)
\right )
,\ a \in {\bold R}^{*}
$$
so that
$$
c^{\prime}(y)^2 e^{-\omega(x)}+
\left ( \frac{\omega^{\prime}(y)}{2} \right )^2
=a^2. 
$$
With respect to this isothermic 
coordinate system,
the mean curvature is given by
$$
H(y)=\frac{1}{8c^{\prime}(y)}
\{4a^2-\omega^{\prime}(y)^2-2\omega^{\prime \prime}(y)
\},\ c^{\prime \prime}(y)=e^{\omega(y)} \omega^{\prime}(y)H(y).
$$
\endexample

\example{Example 4.4}({\bf Spacelike axis and timelike profile curve.})
Let $F:M \longrightarrow {\bold E}^3_1$ be a timelike
surface of revolution with spacelike axis
and timelike profile curve.
Then there exists an isothermic parametrization 
$$
F(x,y)=\frac{1}{a}\left(
e^{\frac{\omega(x)}{2}}\cosh (ay),
e^{\frac{\omega(x)}{2}}\sinh (ay),
c(x)
\right )
,\ a \in {\bold R}^{*}
$$
so that
$$
-c^{\prime}(x)^2 e^{-\omega(x)}+
\left ( \frac{\omega^{\prime}(x)}{2} \right )^2
=a^2. 
$$
With respect to this isothermic coordinate
system,
the mean curvature is given by
$$
H(x)=\frac{1}{8c^{\prime}(x)}
\{4a^2-\omega^{\prime}(x)^2-2\omega^{\prime \prime}(x)
\}
,\ c^{\prime \prime}(x)=-e^{\omega(x)} \omega^{\prime}(x)H(x).
$$
\endexample   

\example{Example 4.5}({\bf Timelike axis.})
Let $F:M \longrightarrow {\bold E}^3_1$ be a timelike
surface of revolution with timelike axis.
Then there exists an isothermic parametrization 
$$
F(x,y)=\frac{1}{a}\left(c(x),
e^{\frac{\omega(x)}{2}}\cos (ay),
e^{\frac{\omega(x)}{2}}\sin (ay),
\right )
,\ a \in {\bold R}^{*}
$$
so that
$$
c^{\prime}(x)^2e^{-\omega(x)}-
\left ( \frac{\omega^{\prime}(x)}{2} \right )^2
=a^2. 
$$
With respect to this isothermic 
coordinate system,
the mean curvature is given by
$$
H(x)=-\frac{1}{8c^{\prime}(x)}
\{2\omega^{\prime \prime}(x)+\omega^{\prime}(x)^2+4a^2
\},\ c^{\prime \prime}(x)=-e^{\omega(x)} \omega^{\prime}(x)H(x).
$$
\endexample   

\example{Example 4.6}({\bf Null axis.})
Let $F:M\longrightarrow {\bold E}^3_1$ be a timelike
surface of revolution with null axis.
Then there exists a null basis 
$\{L_1,L_2,L_3\}$ of ${\bold E}^3_1$ 
and an isothermic parametrization 
$$
F(x,y)=\left (a(x),b(x)-\frac{y^2}{2}a(x),ya(x)
\right )
,\ a \in {\bold R}^{*}
$$
relative to the null basis $\{L_1,L_2,L_3\}$
so that
$$
2a^{\prime}(x)b^{\prime}(x)=-a(x)^2. 
$$
Here a linear null frame 
means a basis of ${\bold E}^3_1$
such that 
$$
\langle L_1,L_1 \rangle=
\langle L_2,L_2 \rangle=0, \ \ 
\langle L_1,L_2 \rangle=1,\ \ 
\langle L_3,L_3 \rangle=1, \ \ 
\langle L_1,L_3 \rangle=
\langle L_2,L_3 \rangle=0.
$$
With respect to this isothermic parametrization,
the mean curvature of $F$ is given by
$$
H=\frac{a^{\prime \prime}(x)a(x)+a^{\prime}(x)^2}
{4a(x)^{2}a^{\prime}(x)}.
$$
\endexample

\proclaim{Proposition 4.7}
For any THIMC surface of revolution with nonconstant
mean curvature, there exists an isothermic 
coordinate system $(x,y)$
such that $H(x)=1/x$ or $H(y)=1/y$. 
\endproclaim

\proclaim{Proposition 4.8}
Let $F:M \longrightarrow {\bold E}^3_1$ be a timelike
surface of revolution with spacelike axis
and Euclidean profile curve
parametrized as in Example {\rm 4.3}
with harmonic inverse mean curvature $1/H=y$
and $a=2$.
Then there exists a real valued function $\phi$
such that
$$
e^{\omega(y)}=\frac{y^2}{4} \left\{
\phi^{\prime}(y)+2\sin \phi(y) \right\}^2, \ \ 
c(y)=-\frac{y^2}{4}\{ \phi^{\prime}(y)^2-4\sin^2 \phi(y)\}.
$$
Furthermore $\phi$ is a solution to the 
third Painlev{\'e} equation of trigonometric form{\rm :}

$$
y\left\{
\phi^{\prime \prime}(y)-2\sin (2\phi(y))  \right \}
+\phi^{\prime}(y)+2\sin \phi(y)=0.
\tag4.2
$$
\endproclaim

\proclaim{Proposition 4.9}
Let $F:M \longrightarrow {\bold E}^3_1$ be a timelike
surface of revolution with spacelike axis
and timelike profile curve
parametrized as in Example {\rm 4.4}
with harmonic inverse mean curvature $1/H=x$
and $a=2$.
Then there exists a real valued function $\phi$
such that
$$
e^{\omega(x)}=\frac{x^2}{4} \left\{
\phi^{\prime}(x)-2\sinh \phi(x) \right\}^2, \ \ 
c(x)=\frac{x^2}{4}\{ \phi^{\prime}(x)^2-4\sinh^2 \phi(x)\}.
$$
Furthermore $\phi$ is a solution to the 
third Painlev{\'e} equation of hyperbolic form{\rm :}

$$
x\left\{
\phi^{\prime \prime}(x)-2\sinh (2\phi(x))  \right \}
+\phi^{\prime}(x)\mp 2\sinh \phi(x)=0.
\tag4.3
$$
\endproclaim

\proclaim{Proposition 4.10}
Let $F:M \longrightarrow {\bold E}^3_1$ be a timelike
surface of revolution with timelike axis
parametrized as in Example {\rm 4.5}
with harmonic inverse mean curvature $1/H=x$
and $a=2$.
Then there exists a real valued function $\phi$
such that
$$
e^{\omega(x)}=\frac{x^2}{4} \left\{
\phi^{\prime}(x)+ 2\cosh \phi(x) \right\}^2, \ \
c(x)=\frac{x^2}{4}\{ \phi^{\prime}(x)^2-4\cosh^2 \phi(x)\}.
$$
Furthermore $\phi$ is a solution to the 
ordinary differential 
equation{\rm :}

$$
x\left \{
\phi^{\prime \prime}(x)-2\sinh (2\phi(x))  \right \}
-\phi^{\prime}(x)-2\cosh \phi(x)=0.
\tag4.4
$$
\endproclaim

\example{Remark 4.11}
The ordinary differential equations
(4.2) and (4.3)
are related to the third
Painlev{\' e} equation.
More precisely let $w=w(x)$ be a solution to
the third Painlev{\'e}
equation:
$$
w^{\prime \prime}-\frac{1}{w}(w^{\prime})^2
+\frac{w^{\prime}}{x}-
\frac{\alpha w^2-\alpha}{x}-\frac{\gamma}{w^{3}}-\frac{\gamma}{w}=0
\tag $P_{\roman{I}\!\roman{I}\!\roman{I}}$
$$ 
with unit modulus, {\it i.e.}, 
$w(x)=e^{\sqrt{-1}\psi(x)}$
for some real valued function $\psi(x)$.
Then ($P_{\roman{I}\!\roman{I}\!\roman{I}}$) is equivalent to
the following ordinary differential equation:
$$
x\left \{
\psi^{\prime \prime}(x)+2\gamma\sin (2\psi(x))  \right \}
+\psi^{\prime}(x)+ 2\alpha\sin \psi(x)=0.
$$
If we choose $\alpha=\gamma=1$ then 
we get (4.2).
 In addition, if we complexified the above
third Painlev{\'e} equation in trigonometric form
and put $\psi=\sqrt{-1}\phi$
then $\phi$ satisfies
$$
x\left \{
\phi^{\prime \prime}(x)+2\gamma\sinh (2\phi(x))  \right \}
+\phi^{\prime}(x)+ 2\alpha\sinh \phi(x)=0.
$$ 
If we choose $\alpha=\mp 1$ and $\gamma=-1$
then we get (4.3).
\endexample
Timelike HIMC surfaces of revolution with null
axis can be classified as follows:

\proclaim{Proposition 4.12}
Let $F:M \longrightarrow {\bold E}^3_1$ be a timelike
surface of revolution with null axis
parametrized as in Example {\rm 4.6}
with harmonic inverse mean curvature 
$1/H=4x$.
Then the function $a(x)$
is a solution to the following
ordinary differential equation{\rm :}

$$
x\left \{
a^{\prime \prime}(x)a(x)+a^{\prime}(x)^{2}
\right \}=a^{2}(x)a^{\prime}(x).
\tag4.5
$$
This ordinary differential equation is explicitly
solved by quadratures. In fact
the solution $a(x)$ is given as follows.
$$
12 \int \frac{a}{2a^3+3a^2+c_1} da=
2\log \vert x \vert +c_2,\ c_1,c_2 \in {\bold R}.
\tag 4.6
$$
\endproclaim

\subhead 4.3
\endsubhead
Next, to study timelike surfaces with imaginary principal
curvatures we shall introduce the notion of
{\it anti isothermic surface}. 

\example{Definition 4.13}
Let $F
: M \rightarrow {\frak M}^3_\nu (c)$ be a timelike
surface. A null coordinate system 
$(u,v)$ is said to be 
{\it anti isothermic} if its
Hopf differentials  
take the following form{\rm :}
$$
Q(u,v)=\frac{1}{2}{\frak q}(u,v)\varrho(u)
,\ \ 
R(u,v)=-\frac{1}{2}{\frak q}(u,v)\sigma(v),
\ \ 
\varrho>0,
\sigma>0.
\tag 4.7
$$
In addition 
$(M,F)$ is said to be {\it anti isothermic} 
if there exists an anti isothermic 
coordinate system
around any point of $M$.
\endexample
Note that anti isothermic property 
implies that $M$ has
imaginary principal curvatures.
In ${\frak M}^3_{1}( c)$, $c\geq0$, 
anti isothermic surfaces have 
non negative Gaussian curvature. (See (2.4).)

The following result plays a fundamental
role in the study of isothermic timelike
surfaces and anti isothermic 
timelike surfaces
in ${\bold E}^3_1$.
We write these alternatives 
together
as $\pm$ isothermic.

\proclaim{Proposition 4.14}
Let $(M,F)$ be a $\pm$ isothermic timelike surface in
${\bold E}^3_{1}$ and
\flushpar
$({\frak D};u,v)$ a simply connected $\pm$
isothermic coordinate region
so that
$$
I=e^\omega dudv,\ \
Q=\frac{1}{2}{\frak q}(u,v)\varrho(u),\
\
R=\pm
\frac{1}{2}{\frak q}(u,v)\sigma(v),\ \ 
\varrho>0, \sigma >0.
$$
Then the formulas{\rm :}
$$
F^{*}_{u}=e^{-\omega}\varrho F_{v},\ 
F^{*}_{v}=\pm e^{-\omega}\sigma F_{u},\ N^{*}=N
\tag{4.8}
$$
define a $\pm$  isothermic timelike immersion
$F^* :{\frak D}\rightarrow {\bold E}^3_1$.
The conformal 
structure of ${\frak D}$ 
induced by 
$F^*$ 
is anti conformal to the original
conformal structure determined by $I$. The
fundamental quantities of
$F^*$
are given as follows{\rm :}
$$
I^{*}=\pm e^{\omega^{*}}dudv=\pm e^{-\omega}\varrho
\sigma dudv,\ \ 
H^*={\frak q}, \ \ Q^*=\varrho H/2,
\ \ R^*=\pm \sigma H/2.
\tag{4.9}
$$
The new immersion $F^*$ is called
the Christoffel transform of $F$ or
dual of $F$. 
\endproclaim

In particular for $\pm$ isothermic THIMC surfaces,
we have the following.

\proclaim{Corollary 4.15}
Every $\pm$ isothermic THIMC surface 
in ${\bold E}^3_1$ is dual to
a $\pm$ timelike Bonnet surface in ${\bold E}^3_1$
and vice versa.
\endproclaim

\head 5. The Hazzidakis equation
\endhead

In this section we shall investigate normal forms
of Gauss equation for THIMC surfaces.
\medskip

\subhead 5.1
Timelike surfaces with $\pm$ holomorphic inverse mean
curvature
\endsubhead
\medskip
\par
Let $F:M\rightarrow {\bold E}^3_1$ be a 
$\pm$ isothermic timelike
surface with $\pm$ holomorphic inverse
mean curvature. 

Without loss of generality we may assume that
$1/H=g(v)$.
Take a $\pm$ isothermic coordinate system $(u,v)$
such that $Q=\varepsilon R={\frak q}(u,v)/2$.
Here $\varepsilon$ denotes the signature $+$ or $-$. 
The Codazzi equations
(C$_0$) become
$$
{\frak q}={\frak q}(u),\ 
e^\omega=
-\frac{\varepsilon g^2{\frak q}_u}{g_v}.
\tag 5.1
$$
Hence we get $\omega_{uv}=0$
and hence $M$ is flat by (2.3).
On the other hand the Gauss equation (G$_0$) implies
$$
e^{2\omega}=\frac{4\varepsilon Q^2}{H^2}
=\frac{\varepsilon {\frak q}^2}{H^2}.
\tag 5.2
$$
Hence $M$ is isothermic.
Moreover (5.1) and (5.2) imply that
$$
g^2{\frak q}^2_u={\frak q}^2g_v^2.
$$
Hence $g{\frak q}_u=\pm {\frak q} g_v$.
Thus we have
$$
g(v)=C_1 e^{\alpha v},\ {\frak q}(u)=C_2 e^{\mu \alpha u},\ \mu=\pm 1,
\ \ 
\alpha \in {\bold R},\
C_1, C_2 \in {\bold R}^{*},\ C_1C_2/\mu<0.
$$
These formulas show that timelike
surfaces with 
$\pm$ holomorphic inverse mean curvature
are flat Bonnet surfaces
with $\pm$holomorphic mean curvature
described in \cite{11, Theorem 3.1}.
In particular the case $\alpha=0$
corresponds to  
timelike CMC cylinders.

\proclaim{Proposition 5.1}
Let $M$ be a timelike  
surface in ${\bold E}^3_1$
with $\pm$ holomorphic inverse mean curvature.
If $M$ is $\pm$ isothermic then
$M$ is a flat isothermic timelike Bonnet surface.
\endproclaim

The notion of $\pm$ isothermic surfaceual can be generalized
to the notion of ``$(\varepsilon,\vartheta)$-i\-so\-ther\-mic surface"
in the following way:

\example{Definition 5.2}
A timelike surface $(M,F)$ is said to be 
$(\varepsilon,\vartheta)$-{\it isothermic}  if there
exists a local null coordinate system
$(u,v)$ around
any point of $M$ such that the Hopf differentials 
$Q$ and $R$ have the following form:
$$
Q(u,v)=\frac{1}{2}({\frak q}(u,v)+\vartheta)\varrho(u),
\ \
R(u,v)=\frac{\varepsilon}{2}({\frak q}(u,v)-\vartheta)\sigma(v),
\ \
\varrho>0,\ \sigma>0.
\tag5.3
$$ 

Here ${\frak q}$ is a real smooth function,
$\varrho$ and $\sigma $ are $\pm$ Lorentz-holomorphic
functions and $\vartheta$ is a real constant.
If $\varepsilon=+$ [resp\. $\varepsilon =-$],
then 
we call
$M$ a $\vartheta$--{\it isothermic} surface
[resp\. an anti $\vartheta$--{\it isothermic} surface].
\endexample

Note that the constant $\vartheta$ has no global meaning,
in fact, $\vartheta$ depends on the choice of $(u,v)$.

\proclaim{Proposition 5.3}
Let $M$ be an $(\varepsilon,\vartheta)$--isothermic
timelike surface with $\vartheta \not=0$. Then
$M$ is $\pm$ isothermic if and only if
$M$ is a timelike Bonnet surface. 
\endproclaim

Proposition 5.1 is generalized as follows:

\proclaim{Theorem 5.4}
Let $M$ be an $(\varepsilon,\vartheta)$--isothermic
timelike surface with Lorentz anti holomorphic
inverse mean curvature $1/H=g(v)$, $\vartheta\not=0$.
Then $M$ is flat and has real distinct principal
curvatures.
\flushpar
{\rm (}{\rm 1}{\rm )} If $M$ is $\vartheta$--isothermic
then $g(v)=C e^{\alpha v},\ {\frak q}(u)=\vartheta \cosh 
(\alpha u+\beta)$,
\medskip
\flushpar
{\rm (}{\rm 2}{\rm )} If $M$ is anti $\vartheta$--isothermic
then $g(v)=C e^{\alpha v},\ {\frak q}(u)=\vartheta \sin 
(\alpha u+\beta),\  C \in {\bold R}^{*},\ 
\alpha,\beta \in {\bold R}$.
\endproclaim

For any $(\varepsilon,\vartheta)$--isothermic
THIMC surface in ${\bold E}^3_1$,
we can consider the {\it dual\/} Bonnet surface in $H^3_1$
or $S^3_1$.

\proclaim{Proposition 5.5}
Let $(M,F)$ be an $(\varepsilon,\vartheta)$-isothermic
timelike surface
in ${\bold E}^3_1$ and
$({\frak D};u,v)$
a simply connected $(\varepsilon,\vartheta)$-isothermic
coordinate region such that 
the Hopf differentials  take the
following forms{\rm :}
$$
Q=\frac{1}{2}({\frak q}(u,v)+\vartheta),\ 
R=\frac{\varepsilon}{2}({\frak q}(u,v)-\vartheta).
$$
Then 
\flushpar {\rm (1)} if $\varepsilon=+$,
there exists a timelike immersion
$$
F^{*}:{\frak D}\longrightarrow
\cases  H^3_1(\frac{1}{|\vartheta|}), & \ \vartheta \not=0, \\
 {\bold E}^3_1, & \ \vartheta=0. 
\endcases
$$
\medskip
\flushpar{\rm (2)}
if $\varepsilon=-$,
there exists a timelike immersion
$$
F^{*}:{\frak D}\longrightarrow
\cases  S^3_1(\frac{1}{|\vartheta|}), & \ \vartheta \not=0, \\
 {\bold E}^3_1, & \ \vartheta=0. 
\endcases
$$
The timelike immersion $F^{*}$ is called a dual surface
of $F$. In particular if F is a THIMC
surface then $F^{*}$ is a timelike Bonnet surface
and vice versa. 
\endproclaim

\example{Remark}
In section 6, we shall prove 
a Lawson correspondence between 
THIMC surfaces in Lorentzian space forms.
Combining the duality in the preceding proposition
and Lawson correspondence, we get a duality between
THIMC surfaces and timelike Bonnet surfaces in $H^3_1$.
\endexample
\subhead 5.2.
Timelike surfaces with non $\pm$ holomorphic inverse mean
curvature
\endsubhead
\medskip

Let $F:M\rightarrow {\bold E}^3_1$ be a THIMC
surface parametrized by a null coordinate 
system $({\bar u},{\bar v})$. 
Since the reciprocal of mean curvature
of $(M,F)$ is harmonic, the mean curvature $H$
can be written as 
$$
\frac{1}{H}=f({\bar u})+g({\bar v}). \tag5.4
$$

Inserting (5.4) in the Codazzi equation
$({\roman C}_0)$ we get
$$
f_{\bar u} R_{\bar u}=g_{\bar v} Q_{\bar v}.
\tag5.5
$$
Inserting this formula into the 
Gauss equation (${\roman G}_0$)
we get
$$
f_{\bar u} \left( \frac
{Q_{{\bar u}{\bar v}}}{Q_{\bar v}} \right )_{{\bar v}}-
Q_{\bar v}
=\frac{f_{\bar u} g_{\bar v}}{(f+g)^2}\left (
2f_{\bar u}-\frac{QR}{R_{\bar u}}
\right ). \tag 5.6
$$
Thanks to (5.5), the equation (5.6) is equivalent
to
$$
g_{\bar v} \left( \frac
{R_{{\bar u}{\bar v}}}{R_{\bar u}} \right )_{{\bar u}}-
R_{\bar u}
=\frac{f_{\bar u} g_{\bar v}}{(f+g)^2}\left (
2g_{\bar v}-\frac{QR}{Q_{\bar v}}
\right ). \tag 5.7
$$

\bigskip

As long as $f_{\bar u} \not=0,\
g_{\bar v}\not=0$, we may assume $\xi:=f({\bar u}),\ \eta:=g({\bar v})$
is a local null coordinate system. With respect to
$(\xi,\eta)$, Gauss-Codazzi equations $({\roman G}_0)$ and
$({\roman C}_0)$
become:
$$
\left( \frac
{Q_{\xi \eta}}{Q_\eta} \right )_{\eta}-
Q_\eta
=\frac{1}{(\xi+\eta)^2}\left (
2-\frac{QR}{R_\xi}
\right ),\ \ 
Q_\eta=R_\xi. \tag 5.8
$$
We should remark that every solution $\{Q,R\}$ to
$$
2-\frac{QR}{R_\xi}=0
\tag5.9
$$
solves (5.8).
Let $\{Q,R\}$ be a solution to (5.9). Then 
by the Codazzi equations $({\roman C}_0)$
and the formula $1/H=\xi+\eta$,
we get
$$
e^{\omega(\xi,\eta)}
=-2(\xi+\eta)^2R_\xi=-
(\xi+\eta)^2
Q(\xi,\eta)
R(\xi,\eta).
$$
Hence the solution $\{Q,R\}$
to $(5.9)$ defines a THIMC
surface if and only if $QR<0$.
Such THIMC surfaces 
have no Euclidean counterparts.
(Compare with Euclidean case \cite{5, p.~203}.)
\smallskip
\par
Hereafter we restrict our attention to 
$(\varepsilon,\vartheta)$--isothermic
THIMC surfaces.
Namely we assume 
$$
Q(\xi,\eta)=\frac{1}{2}
({\frak q}(\xi,\eta)+\vartheta)
\varrho(\xi)
,\ \ 
R(\xi,\eta)=\frac{\varepsilon}{2}
({\frak q}(\xi,\eta)-\vartheta)
\sigma(\eta),\
\ \varrho>0, \sigma>0.
\tag5.10
$$
To adapt our computations to \cite{5} and \cite{10},
and avoid a plethora of 
unnecessary $1/2$'s in the description, 
we shall use the following convention:
$$
q(u,v):=\frac{\varepsilon}{2}{\frak q}(u,v),\ \ \theta:=\frac{1}{2}\vartheta.
$$
And we call $(\xi,\eta)$ simply an
$(\varepsilon,\theta)$--{\it isothermic coordinate system}.

\medskip
Inserting (5.10) to (5.5), we get
$$
\varepsilon \ 
\sigma(\eta)q_{\xi}
(\xi,\eta)
=\varrho(\xi)q_{\eta}
(\xi,\eta).
\tag5.11
$$
Now we introduce a new null coordinate system
$(u,v)$ by
$$
u=\int \varrho(\xi) d\xi,\ \ 
v=\int \sigma (\eta) d\eta.
$$
Then the formula (5.11) implies
that $q$ depends only on 
$t:=\varepsilon u+v$.

We should separate our consideration
to the following two cases:
\bigskip
\flushpar
(1) $
\ \ 2-QR/R_\xi=0,\
\ 
(2)
\ \ 2-QR/R_\xi
\not=0.
$
\bigskip
\subhead 5.3
\endsubhead 
Case 1$:\ 2-QR/R_\xi=0$
\bigskip
In this case, 
the Hopf differentials are given by
$$
Q(\xi,\eta)=\varrho(\xi)
(\varepsilon q(t)+\theta),\ \ 
R(\xi,\eta)=
\sigma(\eta)
(q(t)-\varepsilon \theta).
$$
$$
q(t)=
\cases
-\theta \tanh 
\left (
\theta t/2
\right ),
\ \  \theta \not=0
\\
-2/t,\ \  \theta=0.
\endcases
\tag5.12
$$
 
Inserting (5.12) into 
(C${}_{0}$), we have
$$
e^{\omega(u,v)}
=\cases
\varepsilon \theta^2\left ( 
\xi(u)+
\eta(v)\right )^2
/\cosh^2 (\theta t/2),
\ \ \theta\not=0, \\
-4\varepsilon \left ( 
\xi(u)+
\eta(v)
\right )^2/
t^2,\ \ 
\theta=0.
\endcases
$$
These formulas imply that 
$\varepsilon=+$ for $\theta\not=0$ and
$\varepsilon=-$ for $\theta=0$.

\proclaim{Proposition 5.7}
Let $(M,F)$ be an 
$(\varepsilon,\theta)$--isothermic THIMC
surface in ${\bold E}^3_1$ with 
\flushpar
$(\varepsilon,\theta)$--isothermic coordinate
$(\xi,\eta)$ of the form $(5.10)$ and $(5.11)$.
If $2R_\xi-QR=0$,
\flushpar
then
$\varepsilon=+$ for $\theta\not=0$ and
$\varepsilon=-$ for $\theta=0$.
The fundamental quantities of
$(M,F)$ are given by
$$
Q(u,v)=\frac{\varepsilon q(t)+\theta}
{\varrho (\xi(u))},\ \ \ 
R(u,v)=\frac{q(t)-\varepsilon \theta}
{\sigma(\eta(v))},\ \ \ 
q(t)=
\cases
-\theta \tanh (
\frac{\theta t}{2}), \ \  \theta \not=0
\\
-2/t,\ \  \theta=0,
\endcases
$$
$$
H(u,v)=\frac
{1}{\xi(u)+\eta(v)},\ \ 
I=
\cases
\theta^2(\xi(u)+\eta(v))^2 dudv/
\cosh^{2}(\frac{\theta t}{2}),
\ \ \theta \not=0, \\
4(\xi(u)+\eta(v))^{2}\ dudv/t^2,
\ \  \theta=0.
\endcases
$$
The dual surface of $(M,F)$ is given by
the following formulas
{\rm :}
\smallskip
\flushpar
\rm{(1)} If $\theta \not=0$ then the 
dual surface
$F^{*}$ in $H^3_1(1/(2|\theta|))$
is defined by the data
{\rm :}
$$
e^{\omega^{*}(u,v)}=
\frac{\cosh^{2}(\frac{\theta t}{2})}
{\theta^2(\xi(u)+\eta(v))^2},
$$
$$
Q^{*}(u,v)=R^{*}(u,v)=
\frac{1}{2(\xi(u)+\eta(v))},
\ \ 
H^{*}(u,v)=-2\theta
\tanh (\frac{\theta t}
{2}).
$$
The dual surface 
$F^{*}$ is an isothermic timelike Bonnet surface in 
$H^3_1(1/(2|\theta|))$.
\smallskip
\flushpar
\rm{(2)} If $\theta=0$ then the 
dual surface
$F^*$ in ${\bold E}^3_1$ is defined by
the data{\rm :}
$$
e^{\omega^{*}(u,v)}=
\frac{t^2}{4(\xi(u)+\eta(v))^2},
$$
$$
Q^{*}(u,v)=-R^{*}(u,v)=
\frac{1}{2(\xi(u)+\eta(v))},\ \ 
H^{*}(u,v)
=\frac{-4}{t}.
$$ 
The dual surface $F^*$ is an anti isothermic
timelike Bonnet surface in ${\bold E}^3_1$.
\endproclaim

We call a THIMC surface $(M,F)$ 
{\it generic} if $(M,F)$ does not correspond to
a solution of $2R_\xi-QR=0$. 

\bigskip
\subhead 5.4
\endsubhead
Case 2$:\ 2R_{\xi}-QR\not=0$
\par
In this case, inserting
$$
Q(\xi,\eta)=\varrho(\xi)(\varepsilon
q(t)+\theta),\ \ 
R(\xi,\eta)=\sigma(\eta)(q(t)-\varepsilon
\theta), 
$$
into (5.8) and by the assumption
$2R_{\xi}-QR\not=0$,
we can define the following function
$$
{\Cal S}(t)=
\frac{1}
{\varrho(\xi(u))
\sigma(\eta(v))
(\xi (u)
+\eta (v))^2}.
\tag5.13
$$
\smallskip
The following theorem is proved by much the same way in
\cite{5} and \cite{10}.

\proclaim{Theorem 5.8}
There exist three classes-- $A$, $B$ and $C$-- 
of associated families of generic 
$(\varepsilon,\theta)$--isothermic 
THIMC surfaces in ${\bold E}^3_1$. 
The immersion function of each family is 
given by the Sym formula
$(3.8)$ and $(3.9)$ in Proposition $3.4$, where the data 
$(\omega,Q,R,H)$ in $(3.8)$ are determined by
$$
e^{\omega(u,v)}=
-2\varepsilon q^{\prime}(t)(\xi(u)+\eta(v))^2,
$$
$$
Q(u,v)=\frac{\varepsilon q(t)+\theta}
{\varrho(\xi(u))},\ \
R(u,v)=\frac{q(t)-\varepsilon \theta}
{\sigma(\eta(v))},\ \ 
H(u,v)=\frac{1}{\xi(u)+\eta(v)}.
$$ 
Here $q(t)$ is a solution to 
the generalized Hazzidakis equation{\rm:}
$$
\left(
\frac{q^{\prime \prime}(t)}
{q^{\prime}(t)}
\right )^{\prime}
-q^{\prime}(t)=
{\Cal S}(t)
\left (
2-\frac{q^{2}(t)-\theta^2}{q^{\prime}(t)}
\right),
\ \ -\varepsilon q^{\prime}(t)>0.
\tag$\bigstar_{-\theta^2}^{-\varepsilon}$
$$ 
Here the coefficient function 
${\Cal S}(t)$ in the 
generalized Hazzidakis equation is given by

$$
\vbox{\offinterlineskip
\halign{&\vrule#&\strut\quad#\hfil\quad
\cr
\noalign{\hrule}
height2pt & \omit && 
\omit & \cr
& \hfil Family && \hfil 
Coefficient
& \cr
height2pt & \omit &&
\omit & \cr
\noalign{\hrule}
height2pt & \omit && \omit & \cr
& $A$--family  &&
${\Cal S}(t)=1/\sin^{2}(2t)$ & \cr
\noalign{\hrule}
height2pt & \omit && \omit & \cr
& $B$--family &&
${\Cal S}(t)=1/\sinh^{2}(2t)$ & \cr
\noalign{\hrule}
height2pt & \omit && \omit & \cr
& $C$--family &&
${\Cal S}(t)=1/t^2$ & \cr
height2pt & \omit && \omit & \cr
\noalign{\hrule}
}}
$$

\medskip
\par
Any generic $(\varepsilon,\theta)$--isothermic
THIMC surface belongs to one of these families
$A,B$ or $C$. 
\endproclaim

Via the duality between $\pm 1$--isothermic 
THIMC surfaces in ${\bold E}^3_1$ and 
isothermic timelike Bonnet surfaces in $H^3_1$,
the generalized Hazzidakis equation $(\bigstar^{-}_{-1})$
coincides with that for 
isothermic timelike Bonnet surfaces in $H^3_1$
obtained in \cite{11, Theorem 6.1}.

Moreover the generalized Hazzidakis 
equation $(\bigstar^{-}_{-1})$
coincides with that for 
Bonnet surfaces in hyperbolic $3$-space
$H^3$. 
(See \cite{4, Theorem 3.3.1} and \cite{20}.)
Thus $(\bigstar^{-}_{-1})$ for $A$ or $B$-family
[resp. $C$-family] is solved by 
Painlev{\'e} transcendents $P_{\roman{V}\!\roman{I}}$
[resp. $P_{\roman{V}}$].
See \cite{4, Theorem 3.5.1, 3.5.2}. 
Hence Bonnet surfaces in $H^3$ of non-Willmore type, 
$(\varepsilon,\pm 1)$--isothermic THIMC surfaces 
in ${\bold E}^3_1$
and (generic) timelike Bonnet surfaces 
in $H^3_1$ are 
derived from $P_{\roman V}$ and 
$P_{\roman{V}\!\roman{I}}$.
Note that $(\bigstar^{+}_{\theta^2})$ coincides with
generalized Hazzidakis equation for $\theta$--isothermic 
spacelike HIMC surfaces
in 
${\bold E}^3_1$ 
(and hence spacelike Bonnet surfaces in 
$H^3_1$)
\cite{10}.


\head 6. Timelike HIMC surfaces in ${\frak M}^3_1(c)$
\endhead
\subhead 6.1
\endsubhead
In this section we shall generalize the notion of 
THIMC surface 
in Minkowski $3$-space to that of ${\frak M}^3_1(c)$.

\proclaim{Proposition 6.1}
Let $I[c]$ be a 1-dimensional Riemannian manifold defined by

$$
I[c]=\cases ({\bold R},g[c]\ ) & c=0, \  1, \\
( {\bold R}\setminus \{\pm 1\}, g[c]\ ) & c=-1,
\endcases,\ \ 
g[c] =\frac{dt^2}{(1+ct^2)^2}.
$$

Let $\varphi:M\rightarrow I[c]$ be a
smooth map from a Lorentz surface $M$. Then $\varphi$
is a {\rm (}Lorentzian{\rm )} harmonic map if and only if

$$
\frac{\partial^2\varphi}{\partial u\partial v}-
\frac{2c\varphi}{1+c\varphi^2}
 \frac{\partial \varphi}
{\partial u}\frac{\partial \varphi}
{\partial v}
=0 \tag{6.1}
$$
with respect to any 
{\rm (}and hence in turn all {\rm )} null coordinate
system $(u,v)$.
\endproclaim
The harmonic map equation (6.1) may be considered as
a nonlinear generalization of the
classical linear wave equation $\varphi_{uv}=0$.
As is well known classical linear wave equation
can be solved by the {\it d'Alembert formula}.
The following is regarded as a nonlinear 
d'Alembert formula for (6.1).
\proclaim{Proposition 6.2}
The harmonic map
equation {\rm (6.1)} can be solved as follows{\rm :}

$$
\varphi (u,v)=\cases f(u)+g(v), & c=0, \\
\frac{f(u)+g(v)}{1-cf(u)g(v)} \ {\roman o}{\roman r} 
\ \frac{1-cf(u)g(v)}{f(u)+g(v)},  & c=\pm 1.
\endcases
$$
\endproclaim

The following definition is
 a generalization of that in   
Section 3.

\definition{Definition 6.3}
Let $F:M \rightarrow {\frak M}^3_1(c)$ be a timelike surface.
Then $M$ is said to be a {\it  timelike surface 
with harmonic inverse mean curvature} (THIMC surface)
if $1/H$ is a harmonic
map into $I[c]$. 
\enddefinition  

\medskip
\subhead 6.2
\endsubhead
Hereafter we assume that $M$ is {\it simply connected}.
We denote ${\Cal C}_{H}$ the {\it moduli space} of
conformal immersions of $M$ into ${\frak M}^3_1(c)$
with prescribed mean curvature $H$:
\medskip
\par
${\Cal C}_{H}=
\{ F:M \rightarrow {\frak M}^3_1(c) \ |$ a conformal
timelike immersion 
with mean curvature $H$ $\} / {\frak I}_0(c)$.
\medskip
\flushpar
Here ${\frak I}_0(c)$
is the identity component of the full isometry group
of ${\frak M}^3_1(c)$.
Then we can deduce (by the fundamental theorem of surface theory) that
$$
{\Cal C}_{H} \cong \left \{ (\omega,Q,R) | 
\ {\roman a} \
{\roman s}{\roman o}
{\roman l}{\roman u}{\roman t}{\roman i}
{\roman o}{\roman n}\  {\roman t}{\roman o}\ 
{\roman ( }
{\roman G }_c 
{\roman )} \ 
\roman{and}
\
{\roman (}
{\roman C }_c{\roman ) } \
{\roman w}{\roman i}{\roman t}{\roman h} \
\roman{mean}\
\roman{curvature}\
H \ \right \}. 
$$

\proclaim{Theorem 6.4}{\rm (generalized Lawson correspondences)}
Let $M$ be a simply connected Lorentz surface, f a holomorphic
function and $g$ an anti holomorphic
function on $M$. We define a function $H_c$ by 
$H_c:=(1-cfg)/(f+g)$. 
Then the three moduli
spaces 
$
{\Cal C}_{H_0}, \ {\Cal C}_{H_1}, 
\ {\Cal C}_{H_{-1}}$
are
mutually isomorphic.

\endproclaim
\demo{\bf Proof}
Let $(\omega, Q,R, H_c)$ 
be a solution of 
$({\roman G}_c)$ and $({\roman C}_c)$ for 
$ c=\pm 1 $. 
Then 
$({\tilde \omega},{\tilde Q},{\tilde R},H_0)$
defined by 
$$
e^{{\tilde \omega}}:=(1+cf^2)(1+cg^2)\ e^{\omega},
\ \ {\tilde Q}:=(1+cf^2)Q,
\ \ {\tilde R}:=(1+cg^2)R
$$
is a solution to $({\roman G}_c)$ and 
$({\roman C}_c)$. 
Note that in case $c=-1$, the function $(1+cf^2)(1+cg^2)$
is positive if and only if $H_{-1}^2>1$.
\qed
\enddemo

Theorem 6.4 may be considered as a generalization of the so-called 
{\it Lawson correspondences} for timelike CMC surfaces.

\example{Remark}
In Riemannian case,
(non CMC) HIMC surfaces in $H^3$ have Lawson correspondents
if and only if $H^2>1$. 
On the other hand,
in spacelike case, 
(non CMC) spacelike HIMC surfaces
in $S^3_1$ have Lawson correspondents
if and only if $H^2>1$. 
See \cite{7}, \cite{10}.
\endexample
\subhead 6.3
\endsubhead
Using the Lawson correspondences described above,
we can give immersion formulas for THIMC surfaces in ${\frak M}^3_1(c)$,
$c=\pm 1$. 
\par
Before describing immersion formulas for THIMC surfaces, 
we point out the following invariance of (6.1):
\smallskip
Let $\varphi$ be a solution to (6.1) of the form:
$$
\varphi(u,v)=\frac{f(u)+g(v)}{1-cf(u)g(v)}.
$$

Then the replacements:
$$
f \longmapsto  2\tau f,\  g \longmapsto 
2\tau g,\ \ \tau \in {\bold R}^{*}
$$
produce a new solution to (6.1).
More precisely, the function $\varphi[\tau]$ defined by
$$
\varphi[\tau](u,v)=\frac{2\tau(f(u)+g(v))}{1-4c\tau^{2}f(u)g(v)}
$$
is still a solution to (6.1).
 
Let $\Phi_{\lambda}$ be a solution of the zero curvature
equations (3.8) with variable
spectral parameter $\lambda$.
To describe immersion formulas we shall use
the following notational convention.
$$
\Phi[\tau]:=\Phi_{\lambda},\ \lambda=
(1-2\tau g)/
(1+2\tau f),\ \  \tau \in {\bold R}.
$$ 
Since
the zero curvature equation (3.8) is completely integrable,
(3.8) has also solutions for all 
$\tau \in {\bold C}$.

Direct computations similar to those in \cite{1}, \cite{7}
and \cite{10}
show the following.
\proclaim{Theorem 6.5}{\rm (Immersion formulas)}\flushpar

Let $\Phi[\tau]: M \times {\bold C} \rightarrow 
G^{\bold C}$ be a complexified
solution to {\rm (}{\rm 3}.{\rm 8}{\rm )}.
Then the followings hold.
\flushpar
\smallskip
{\rm (}$c=0${\rm )} For every $\tau \in {\bold R},$
$$
F^{(0)}(\tau):=-\frac{\partial}{\partial \tau}\Phi[\tau]
\cdot \Phi[\tau]^{-1},\ \tau \in {\bold R}
$$ 
describes a THIMC surface in 
${\bold E}^3_1$ given in {\rm Proposition 3.4}.
\smallskip
\flushpar
\smallskip
{\rm (}$c=-1${\rm )} For any 
$\tau \in
{\bold R}^{*}$, 
$$
F^{(-1)}(\tau):=p_H( \Phi[\tau], \Phi[-\tau])
$$
is a THIMC surface in $H^3_1$
with 
unit normal vector field
$$
N=-\mu_{H}(\Phi[\tau], \Phi[-\tau]){\bold k}^{\prime}.
$$
\smallskip
\flushpar
\smallskip
{\rm (}$c=1${\rm )} Let $\Phi[\sqrt{-1}\tau],
\ \tau \in {\bold R}$ be a complexified solution to 
{\rm (3.8)}. 
\flushpar
Then for every $\tau \in {\bold R}$
$$
F^{(1)}(\tau):=p_S(\Phi[\sqrt{-1}\tau])
$$
is a THIMC surface in $S^3_1$
with unit normal vector field
$$
N=-\mu_{S}(\Phi[\sqrt{-1}\tau]) {\bold j}^{\prime}.
$$
The first fundamental form of $F^{(c )},\ c=\pm 1$ is
$$
I^{( c)}(\tau)=\frac{4\tau^2 e^\omega}{(1+4c\tau^2f^2)(1+4c\tau^2g^2)}.
$$
The mean curvature of $F^{(c)}(\tau),\  c=\pm1$ is given by
$$
H=\frac{1-4c\tau^2 fg}{2\tau(f+g)}.
$$
Moreover the mean curvature of $F^{(-1)}(t)$ satisfies
$H^2>1$.

In particular, for $c=\pm1$,
$F^{(\pm1)}(1/2)$
is
the Lawson correspondent
of $F=F^{(0)}$.
The conformal deformations of THIMC surfaces in
${\frak M}^3_1(c)$ preserve
$K/(H^2+c)$.
\endproclaim

\example{Remark} 
The conformal deformation of HIMC surfaces in Riemannian space forms
[resp. spacelike HIMC surfaces in Lorentzian space forms]
preserves $K/(H^2+c)$ [resp. $K/(H^2-c)$]
Note that in case $c=0$, the constancy of $K/(H^2-c)$ is equivalent
to the constancy of the ratio of principal curvatures.
\endexample
Computing the Gaussian curvature or $K/(H^2+c)$,
we have the following theorem.

\proclaim{Theorem 6.7}
Let $(M,F)$ be an 
$(\varepsilon,\theta)$--isothermic
THIMC surface in ${\frak M}^3_1(c )$.
\flushpar {\rm(}{\rm 1}{\rm )} If
$K$ is constant then $K=0$ or $c$.
\smallskip
\flushpar
{\rm(}{\rm 2}{\rm )} 
If $K/(H^2+c)$ is constant then
$(M,F)$ is a flat timelike Bonnet surface.
\endproclaim

As an application of Lawson 
correspondence as above, one can classify
$\pm$ isothermic flat 
timelike Bonnet surfaces in Lorentzian space forms.
In fact, since
the Lawson correspondence preserves 
$\pm$ isothermic property or 
flatness, we have obtained the following (\cite{11, Theorem 6.2}).

\proclaim{Theorem 6.8}
Flat simply connected $\pm$ isothermic timelike Bonnet surfaces in 
\flushpar
one
Lorentzian $3$-space form correspond to those in another 
Lorentzian $3$-space form.
\endproclaim  
\smallskip
Moreover in \cite{11}, timelike Bonnet surfaces
in ${\frak M}^3_1 (c)$
with constant Gaussian curvature
are classified.  
\smallskip

\subhead 6.4
\endsubhead
Finally we consider THIMC surfaces in $H^3_1$ with mean curvature $H^2<1$.
To investigate such surfaces, we use the following invariance of
(6.1) with $c=\pm1$.
Let $H$ be a solution of (6.1) of the form:
$$
H(u,v)=\frac{f(u)+g(v)}{1-cf(u)g(v)}.
$$
Then 
for any $\tau \in {\bold R}^*$, 
the replacements:
$$
f \longmapsto \tau f,\ \ g \longmapsto \tau^{-1}g
$$
produce a new solution of (6.1).
Namely the function
$H[\tau]$ defined by
$$
H[\tau](u,v)=\frac{\tau f(u)+\tau^{-1}g(v)}{1-cf(u)g(v)}
$$
is also a solution of (6.1).
Based on this deformation, we define two auxiliary
functions ({\it variable spectral parameters}):
$$
\lambda(u,\tau):=\frac{\tau(1-cf(u)^2)}{\tau^2-cf(u)^2},\ \
\nu(v,\tau):=\frac{\tau(1-cg(v)^2)}{\tau^2-cg(v)^2}.
$$
Then we have the following. 

\proclaim{Theorem 6.9}
Let $\Psi[\tau]$ be a solution to
$$
\frac{\partial}{\partial u}\Psi[\tau]=
\Psi[\tau]U[\tau],\ \
\frac{\partial}{\partial v}\Psi[\tau]=
\Psi[\tau]V[\tau],
\tag6.2
$$
$$
U[\tau]=\pmatrix
-\frac{1}{4}\omega_u & -Qe^{-\omega/2} \\
\frac{1}{2}(H[\tau]+c)\lambda e^{\omega/2} &
\frac{1}{4}\omega_u
\endpmatrix,\ \ 
V[\tau]=\pmatrix
\frac{1}{4}\omega_v 
& 
-\frac{1}{2}(H[\tau]-c)\nu e^{\omega/2} \\
Re^{-\omega/2} &
-\frac{1}{4}\omega_v
\endpmatrix.
$$
Then for any $\tau \in {\bold R}^{*}$, 
$$
F^{(-1)}[\tau](u,v):=p_{H}(\Psi[\tau],\Psi[-\tau])
$$
is a THIMC surface in $H^3_1$ with unit normal vector field
$N=-\mu_{H}(\Psi[\tau],\Psi[-\tau]){\bold k}^{\prime}$
and mean curvature $H[\tau]$.
The first fundamental form of $F[\tau]$ is given by
$$
I^{(-1)}[\tau]=\frac{\tau^2(1-cf(u)^2)(1-cg(v)^2)e^\omega}
{(\tau^2-cf(u)^2)(\tau^2-cg(v)^2)}dudv.
$$
\endproclaim

Since the Lax equation (6.2) with two variable spectral parameters
$\lambda$ and $\nu$ is completely integrable, (6.2) has solutions for all
$\tau \in {\bold C}$.
Such complexified solutions $\Psi[\tau]$ to (6.2) describe another kind of surfaces in $S^3_1$. 

\proclaim{Theorem 6.10}
Let $\Psi[\tau]:M\times {\bold C} \to G^{\bold C}$ be a 
complexified solution to $(6.2)$. 
Then for any 
$\tau \in {\bold R}^{*}$, 
$$
F^{(1)}[\tau](u,v)=p_{S}(\Psi[\sqrt{-1}\tau])
$$
is a timelike surface in $S^3_1$ with unit normal vector field
$$
N=\mu_{S}(\Psi[\sqrt{-1}\tau]){\bold j}^{\prime}
$$
and mean curvature 
$$
\frac{\tau f-\tau^{-1}g}{1-cfg}.
$$
The first fundamental form of $F^{(1)}[\tau]$ is given by
$$
I^{(1)}[\tau]=\frac{\tau^2(1-cf(u)^2)(1-cg(v)^2)e^\omega}
{(\tau^2+cf(u)^2)(\tau^2+cg(v)^2)}dudv.
$$
\endproclaim
The inverse mean curvature 
of $F^{(1)}[\tau]$ 
in Theorem 6.10 is a harmonic map into $I[-1]$.

\enddocument